\documentclass[a4paper,11pt]{article}

\usepackage{amsmath,amssymb,amsthm}
\usepackage[dvips]{graphicx,color}

\numberwithin{equation}{section}

\newtheorem{thm}{Theorem}[section]

\newtheorem{lem}[thm]{Lemma}
\newtheorem{axthm}{Theorem A}
\newtheorem{uthm}{Theorem B}
\newtheorem{rem}[thm]{Remark}

\allowdisplaybreaks

\begin{document}
\title{On the Best Constant in the Moser-Onofri-Aubin Inequality}
\author{Nassif Ghoussoub\,$^1$ and Chang-Shou Lin\,$^2$}
\date{\small\it $^1$Department of Mathematics, University of British Columbia,\\ Vancouver, BC V6T1Z2, Canada\vspace{.3cm}\\
$^2$Department of Mathematics, Taida Institute for Mathematical Sciences, \\National Taiwan University, Taipei, 106, Taiwan}
\maketitle
\begin{abstract}
Let $S^2$ be the 2-dimensional
unit sphere and let $J_\alpha $ denote the nonlinear functional on the  Sobolev space
$H^{1,2}(S^2)$ defined by
$$
J_\alpha(u) =  \frac{\alpha}{4}\int_{S^2}|\nabla u|^2\, d\omega
+  \int_{S^2} u\, d\omega
-\ln \int_{S^2} e^{u}\, d\omega,
$$
where $d\omega$ denotes Lebesgue measure on  $S^2$, normalized so that
$\int_{S^2} d\omega = 1$.
Onofri had established that $J_\alpha$
  is non-negative on $H^1(S^2)$ provided $\alpha \geq 1$. In this note, we show that if $J_\alpha$ is restricted to those $u\in H^1(S^2)$ that satisfy the Aubin condition:
\begin{equation*}
\int_{S^2}e^u\,x_j\, dw=0\quad\text{for all }1\leq j\leq 3,
\end{equation*}
then the same inequality continues to hold (i.e., $J_\alpha (u)\geq0$) whenever  $\alpha \geq \frac{2}{3}-\epsilon_0$ for some $\epsilon_0>0$.
The question of  Chang-Yang on  whether this remains true for all
$\alpha \geq \frac{1}{2}$ remains open.
\end{abstract}

\section{Introduction}
\hspace{.5cm}
Let $S^2$ be the 2-dimensional unit sphere with the standard metric $g$ and the corresponding volume form $d\omega$ normalized so that  $\int_{S^2} d\omega = 1$. For $\alpha >0$, we consider the following nonlinear functional on  the  Sobolev space
$H^{1,2}(S^2)$:
$$
J_\alpha(u) =  \frac{\alpha}{4}\int_{S^2}|\nabla u|^2\, d\omega
+  \int_{S^2} u\, d\omega
-\ln \int_{S^2} e^{u}\, d\omega.
$$
 The classical Moser-Trudinger inequality \cite{M} yields that  $J_\alpha$ is bounded from below in $H^1(S^2)$ if and only if $\alpha \geq 1$. In \cite{O}, Onofri proved that the infimum is actually equal to zero for  $\alpha =1$, by using the conformal invariance of $J_1$ to show that
\begin{equation}
\inf_{u\in {\mathcal M}}J_1(u)=\inf_{u\in H^1(S^2)}J_1(u)=0,
\end{equation}
where ${\mathcal M}$ is the submanifold of $H^1(S^2)$ defined by
\begin{equation}\label{cond1}
{\mathcal M}:=\left\{u\in H^1(S^2);\, \int_{S^2}e^u{\bf x} \, dw=0 \right\},
\end{equation}
with ${\bf x}= (x_1, x_2, x_3)\in S^2$, on which the infimum of $J_1$ is attained.  Other proofs were also given by Osgood-Phillips-Sarnak \cite{OPS} and by Hong \cite{H}.

Prior to that, Aubin \cite{A} had shown that by restricting the functional $J_\alpha$ to ${\mathcal M}$, it is then again bounded below by ---a necessarily
non-positive--- constant $C_\alpha$,
 for any $\alpha \geq \frac{1}{ 2}$. In their work on Nirenberg's prescribing Gaussian
curvature problem on $S^2$, Chang and Yang \cite{CY, CY2} showed that $C_\alpha$ can be taken
to be equal to $0$ for $\alpha \geq 1-\epsilon_0$ for some small $\epsilon_0$. This led them to the
following\\

\noindent {\bf Conjecture 1:}
 \hbox{If  $\alpha \geq \frac{1}{ 2}$\, then \,$
\inf\limits_{u\in {\mathcal M}} J_\alpha (u)  = 0$.}\\

Note that this fails  if  $\alpha < \frac{1}{2}$, since the functional $J_\alpha$ is then unbounded from
below (see \cite{FFGG}).

\noindent
In this note, we want to give a partial answer to this question by showing that this is indeed the case for $\alpha \geq \frac{2}{ 3}$ and slightly below that.

As mentioned above, Aubin had proved that for all $\alpha \geq \frac{1}{ 2}$, the functional $J_\alpha$ is  coercive on ${\mathcal M}$, and that it attains its infimum  on some function $u\in{\mathcal M}$. Accounting for the Lagrange multipliers,  and setting $\rho=\frac{1}{\alpha}$, the Euler-Lagrangian equation for $u$ is then
\begin{equation*}
\Delta u+2\rho\left(
\frac{e^u}{\int_{S^2}e^u\, dw}-1
\right)=\sum_{j=1}^3\alpha_j x_j e^u\quad\text{on }S^2.
\end{equation*}
In \cite{CY2}, Chang and Yang proved however  that $\alpha_j, j=1,2,3$ necessarily vanish.
 Thus $u$ satisfies -- up to an additive constant -- the following equation:
\begin{equation}\label{eq1}
\Delta u+2\rho(e^u-1)=0\quad\text{on }S^2.
\end{equation}
Conjecture (1) is therefore equivalent to the question whether if $1< \rho \leq 2$, then $u\equiv 0$ is the only solution of \eqref{eq1}.

Here is the  main result of this note.
\begin{thm}\label{thm1}
If $1< \rho\leq\frac{3}{2}$ and $u$ is a solution of \eqref{eq1}, then $u\equiv0\text{ on }S^2$.
\end{thm}
\noindent
This clearly gives a positive answer to the question of Chang and Yang for $\alpha \geq \frac{2}{3}$.

\section{The axially symmetric case}

The proof of Theorem \ref{thm1} relies on the fact that the conjecture has been shown to be true in the axially symmetric case. In other words, the following result holds.
\begin{axthm}
Let $u$ be a solution of \eqref{eq1} with $1<\rho\leq2$. If $u$ is axially symmetric, then $u\equiv0$ on $S^2$.
\end{axthm}
\noindent
Theorem (A) was first established by Feldman, Froese, Ghoussoub and Gui \cite{FFGG} for $1<\rho\leq\frac{25}{16}$. It was eventually proved for all $1<\rho\leq2$ by Gui and Wei \cite{GW}, and  independently by Lin \cite{L}. Note that this means that the following one-dimensional inequality holds:
\begin{equation}  {1\over
2}\int_{-1}^1 (1-x^2)|g'(x)|^2\ dx + 2 \int_{-1}^1 g(x)\ dx -
2\ln{{1}\over{2}}\int_{-1}^1 e^{2g(x)} dx \geq 0,  
\end{equation}
for every function $g$ on $(-1,1)$
satisfying $\int_{-1}^1(1-x^2)|g'(x)|^2 dx <\infty$ and $
\int_{-1}^1 e^{2g(x)} x dx = 0$. \hfill $\square$.\\

We now give a sketch of the proof of Theorem A that connects the conjecture of Chang-Yang to an equally interesting  Liouville type theorem on $R^2$. For that, we let $\Pi$ denote the stereographic projection  $S^2\rightarrow\Bbb R^2$ with respect to the North pole $N=(0,0,1)$:
\begin{equation*}
\Pi(x):=\left(\frac{x_1}{1-x_3}, \frac{x_2}{1-x_3}\right).
\end{equation*}
Suppose $u$ is a solution of \eqref{eq1}, and set
\begin{equation*}
\tilde u(y):=u(\Pi^{-1}(y))\quad\text{for }y\in\Bbb R^2.
\end{equation*}
Then $\tilde u$ satisfies
\begin{equation*}
\Delta \tilde u+8\pi\rho \,J(y)\left(e^{\tilde u}-\frac{1}{4\pi}\right)=0\quad\text{in }\Bbb R^2,
\end{equation*}
where $J(y):=\left(\frac{2}{1+|y|^2}\right)^2$ is the Jacobian of $\Pi$. By letting
\begin{equation}\label{def2}
v(y):=\tilde u(y)+\rho\log\left((1+|y|^2)^{-2}\right)+\log(32\pi\rho)\quad\text{for }y\in\Bbb R^2,
\end{equation}
we have that $v$ satisfies
\begin{equation}\label{eq2}
\Delta v+(1+|y|^2)^l e^v=0\quad\text{in }\Bbb R^2,
\end{equation}
where $l=2(\rho-1)$.\\

Note that by using \eqref{def2} with $u\equiv 0$, equation \eqref{eq2} always has a special axially symmetric solution, namely
\begin{equation}
v^*(y)=-2\rho\log(1+|y|^2)+\log(32\pi\rho)\quad\text{for }y\in\Bbb R^2,
\end{equation}
where again $l=2(\rho-1)$. Moreover,  The Pohozaev idendity yields that for any solution $v$   of \eqref{eq2} we have
\begin{equation}
4<\beta_l (v)<4(1+l),
\end{equation}
where
\begin{equation*}
\beta_l (v):=\frac{1}{2\pi}\int_{\Bbb R^2}(1+|y|^2)^le^{v}dy.
\end{equation*}
An open question that would clearly imply the conjecture of Chang and Yang is the following:\\

\noindent {\bf Conjecture 2:} Is $v^*$ the only solution of \eqref{eq2} whenever $l>0$?\\

Note that it is indeed the case if $\ell <0$ (i.e., $\rho<1$ and $\alpha >1$), since then we can employ the method of moving planes to show that $v(y)$ is radially symmetric with respect to the origin, and then conclude that $u(x)$ is axially symmetric with any line passing through the origin.
Thus $u(x)$ must be a constant function on $S^2$. Equation \eqref{eq1} then yields $u=0$, which implies $J_\alpha \geq 0$ on $\mathcal M$. By passing to the limit as $\alpha\rightarrow1$, we recover the Onofri inequality.\\

When $l>0$ (i.e., $\rho>1$ and $\alpha \leq 1$), the method of moving planes fails and it is still an open problem whether any solution of \eqref{eq2} is equal to $v^*$ or not. The following uniqueness theorem reduces however  the problem to whether any solution of \eqref{eq2} is radially symmetric.

\begin{uthm}\label{thma}
Suppose $l>0$ and $v_i(y)=v_i(|y|), i=1,2,$ are two solutions of \eqref{eq2} satisfying
\begin{equation}
\beta_l (v_1)=\beta_l (v_2).
\end{equation}
Then   $v_1=v_2$ under one of the following conditions:\\
(i) $l\leq1$,

or\\
(ii) $l>1$ and
$2l<\beta_l(v_i)<2(2+l)$ for $ i=1,2$.
\end{uthm}
\noindent

In order to show how Theorem B implies Theorem A, we suppose $u$ is a solution of \eqref{eq1} that  is axially symmetric with respect to some direction. By rotating, the direction can be assumed to be $(0,0,1)$. By
using the stereographic projection as above, and setting $v$ as in  \eqref{def2}, we have
\begin{equation}\label{order}
\begin{cases}
&v(y)=-4\rho\log|y|+O(1),\\
&\frac{1}{2\pi}\int_{\Bbb R^2}(1+|y|^2)^le^vdy=4\rho=4+2l.
\end{cases}
\end{equation}
If $l\leq1$, i.e., $\rho\leq\frac{3}{2}$, then $v=v^*$ by (i) of Theorem B, and then $u\equiv0$. If $l>1$, then by noting that
\begin{equation*}
2l<4\rho=4+2l=\beta_l(v)<4+4l,
\end{equation*}
we deduce that $v=v^*$ by (ii) of Theorem B, which again means that $u\equiv 0$.

\section{Proof of the main theorem}

We shall prove Theorem \ref{thm1} by showing that if $\rho \leq \frac{3}{2}$, then any solution of (\ref{eq1}) is necessarily axially symmetric. We can then conclude by using Theorem A. \\

We shall need the following lemma.
\begin{lem}\label{lem1}
Let $\Omega$ be a simply connected domain in $\Bbb R^2$, and suppose $g\in C^2(\Omega)$    satisfies
\begin{equation*}
\begin{cases}
&\Delta g+e^g>0\quad\text{in }\Omega\quad\text{and}\\
&\int_\Omega e^gdy\leq 8\pi.
\end{cases}
\end{equation*}
Consider an open set $\omega\subset\Omega$ such that $\lambda_{1,g}(\omega)\leq0$, where $\lambda_{1,g}(\omega)$ is  the first eigenvalue of the operator $\Delta +e^g$ on $H^1_0(\omega)$. Then,  we necessarily have  that
\begin{equation}\label{bb}
\int_\omega e^gdy>4\pi.
\end{equation}
\end{lem}
\noindent
Lemma \ref{lem1} was first proved in \cite{B} by using the classical Bol inequality. The strict inequality of \eqref{bb} is due to the fact that $\Delta g+e^g>0$ in $\Omega$. See \cite{BL} and references therein.\\

Now we are in the position to prove the main theorem.\\

\noindent
{\bf Proof of Theorem \ref{thm1}.} Suppose $u(x)$ is a solution of \eqref{eq1}. Let $\xi_0$ be a critical point of $u$. Without loss of generality, we may assume $\xi_0=(0,0,-1)$. By using the stereographic projection $\Pi$ as before and letting
\begin{equation*}
v(y):=u(\Pi^{-1}(x))-2\rho\log(1+|y|^2)+\log(32\pi\rho),
\end{equation*}
$v$ satisfies \eqref{eq2} and
\begin{equation}\label{saiyou}
\nabla v(0)=0.
\end{equation}
Set
\begin{equation*}
\varphi(y):=y_2\frac{\partial v}{\partial y_1}-y_1\frac{\partial v}{\partial y_2}.
\end{equation*}
Then $\varphi$ satisfies
\begin{equation}
\Delta \varphi+(1+|y|^2)^le^v\varphi=0\quad\text{in }\Bbb R^2.
\end{equation}
If $\varphi\not\equiv0$, then by \eqref{saiyou},
\begin{equation*}
\varphi(y)=Q(y)+\text{higher order terms}\quad\text{for }|y|\ll1,
\end{equation*}
where $Q(y)$ is a quadratic polynomial of degree $m$ with $m\geq2$, that is also a harmonic function, i.e., $\Delta Q=0$. Thus, the nodal line $\{y\,|\, \varphi(y)=0\}$ divides a small neighborhood of the origin into at least four regions.  Globally, $\Bbb R^2$
is therefore divided by the nodal line $\{y\,|\,\varphi(y)=0\}$ into at least $3$ regions,  i.e.,
\begin{equation*}\displaystyle
\Bbb R^2\setminus\{y\,|\,\varphi(y)=0\}=\bigcup_{j=1}^3\Omega_j.
\end{equation*}
In each component $\Omega_j$, the first eigenvalue of $\Delta+(1+|y|^2)^le^v$ being equal to $0$. Let now
\begin{equation*}
g:=\log\left((1+|y|^2)^le^v\right).
\end{equation*}
By noting that
\begin{equation*}
\Delta g+e^g>0\quad\text{in }\Bbb R^2,
\end{equation*}
Lemma \ref{lem1} then  implies that for each $j=1, 2, 3$,
\begin{equation*}
\int_{\Omega_j}e^gdy=\int_{\Omega_j}(1+|y|^2)^le^vdy>4\pi.
\end{equation*}
It follows that
\begin{equation*}
8\pi\rho=\sum_{j=1}^3\int_{\Omega_j}(1+|y|^2)^le^vdy>12\pi,
\end{equation*}
which is a contradiction if we had assumed that $\rho\leq \frac{3}{2}$. Thus we have $\varphi(y)=0$, i.e., $v(y)$ is axially symmetric. By Theorem A, we can conclude $u\equiv0$.
\phantom{}\hfill{$\square$}

\begin{rem}\rm
If we further assume that the antipodal of $\xi_0$ is also a critical point of $u$, then $\Bbb R^2\setminus\{y\,|\,\varphi(y)=0\}=\displaystyle\bigcup_{j=1}^m\Omega_j$, where $m\geq4$.
Lemma \ref{lem1} then yields
\begin{equation*}
8\pi\rho=\int_{\Bbb R^2}(1+|y|^2)^le^vdy\geq\sum_{j=1}^m\int_{\Omega_j}(1+|y|^2)^le^vdy>4m\pi\geq16\pi,
\end{equation*}
which is a contradiction whenever $\rho \leq 2$. By Theorem A, we have again that $u\equiv 0$.\\

For example, if $u$ is even on $S^2$ (i.e., $u(z)=u(-z)$ for all $z\in S^2$), then the main theorem holds for $\rho\leq2$.
\end{rem}

\begin{rem} \rm One can actually show that Conjecture 1 holds for $\rho \leq \frac{3}{2}+\epsilon_0$ for some $\epsilon_0>0$.
Indeed, it suffices to show that for $\alpha$ smaller but close to $\frac{2}{3}$, the functional $J_\alpha$ is non-negative. 
Assuming not, then there  exists a sequence of $\{\alpha_k\}_k$  such that $ {1 \over 2} < \alpha_k <{
2\over 3}$, $\lim_k \alpha_k ={2 \over 3} $  and $\inf_{\mathcal M}
J_{\alpha_k}(u) < 0$. Since $J_\alpha$ is coercive for each $\alpha >\frac{1}{2}$, a standard compactness argument yields the existence of a minimizer $u_k \in {\mathcal M}$ for $J_{\alpha_k} $ such that $u_k(0)=0$. Moreover, $\| u_k\|_{H^1} < C$ for
some positive constant independent of $k$. Modulo extracting a
subsequence,  $u_k$ then converges weakly to some $u_0$  in ${\mathcal M}$ as $k \to
\infty$, and $u_0$ is necessarily a minimizer for $I _{{ 2 \over 3}} $ in ${\mathcal M}$ that satisfies $u_0(0) =0$. By our main result, $u_0 \equiv 0$. Now, we
claim that $u_k$ actually converges strongly in $H^1$ to $ u_0 \equiv 0$ . This
is because -- as argued by Chang and Yang --  the Euler-Lagrange equations are then
$$
\frac{\alpha_k}{2} \Delta u_k - 1 +{ 1 \over \lambda_k }e^{u_k}=0
\eqno(3.20)
$$
where $\lambda_k =  \int_{S^2}e^{u_k}dx <C$ for some positive constant $C$.
Multiplying (3.20) by $u_k$ and integrating over $S^2$, we obtain
$$
\frac{\alpha_k}{2} \int_{S^2}|\nabla u_k|^2\,dw  +\int_{S^2}
u_k(x)\ dw = { 1 \over \lambda_k }  \int_{S^2}e^{u_k(x)} u_k(x)\
dw. \eqno(3.21)
$$
Applying Onofri's inequality for $u_k$ and using that $\|u_k\|_{H^1} <C$,  we get
that  $\int_{S^2} e^{2u_k}\ dw$ is also uniformly bounded. This combined with
H\"older's inequality and the fact that $u_k$ converges strongly to $0$ in
$L^2$ yields that $\int_{S^2}e^{u_k} u_k\ dw \to 0$. Use now
(3.21) to conclude that $\|u_k\|_{H^1}  \to 0$ as $k \to \infty$.

Now, write  $u = v + o(||u||)$ for $||u||$ small, where $v$ belongs to the
tangent space of the submanifold ${\mathcal M}$ at $u_0 \equiv 0$ in
$H^1(S^2)$. It is easy to see that $\int_{S^2} v {\bf x}\ dw =0$. We can calculate
the second variation of $J_{\alpha}$ in  ${\mathcal M}$ at $u_0 \equiv 0$ and
get the following estimate around $0$
$$
J_{\alpha} (u)= \alpha \int_{S^2} |\nabla v|^2\,dw
    -2\int_{S^2} |v|^2\,dw + o(||u||^2). 
$$

Note that the eigenvalues of the Laplacian on $S^2$ corresponding to the eigenspace generated by $x_1, x_2, x_3$ are $\lambda_2=\lambda_3=\lambda_4=2$, while  $\lambda_5=6$. Since  $v$ is orthogonal to ${\bf x}$, we have
$$
\int_{S^2} |\nabla v|^2\,dw \ge 6 \int_{S^2} |v|^2\,dw
$$
and therefore
$$ J_{\alpha} (u) \ge (\alpha- {1 \over 3}) ||u||^2 +o(||u||^2).
$$
Taking
$\alpha =\alpha_k$ and $u=u_k$ for $k$ large enough, we get that
$ J_{\alpha_k}(u_k) \ge 0$, which clearly contradicts our initial assumption on $ u_k$.
\end{rem}

\noindent
{\bf Concluding remarks.} (i) The question whether $J_\alpha (u)\geq0$ for $\frac{1}{2}\leq \alpha < \frac{2}{3}$ under the condition \eqref{cond1} is still open. However, in \cite{L2}, it was proved that
there is a constant $C\geq0$ such that for any solution $u$ of \eqref{eq1} with $1<\rho\leq2$ (i.e. $\frac{1}{2}\leq \alpha < 1$), we have
\begin{equation*}
|u(x)|\leq C\quad\text{for all }x\in S^2.
\end{equation*}
\noindent
(ii) Recently, Liouville type equations with singular data have attracted a lot of attentions in the research area of nonlinear partial differential equations, because
it is closely related to vortex condensates appeared in many physics models. One of difficult subjects in this area is to understand bubbling phenomenons arised from solutions of these equations. For the past
twenty years, there have been many works devoted to this direction. Among bubbling phenomenons, the most delicate is the situation when more than one vortex are collapsed into one single point. The equation \eqref{eq2} is one of model equations which can allow us to accurately describe bubbling behavior during those collapses. See
\cite{BLT} and \cite{DET} for related details. Thus, understanding the structure of solutions to the equation \eqref{eq2} is fundamentally important. As mentioned above,  it is conjectured that for $l\leq2$, all solutions of \eqref{eq2} must be radially symmetric. This remains an open question, although a partial answer  has been given recently in \cite{BLT}.

\end{document}